\documentclass[11pt]{article}
\usepackage{amsfonts}
\usepackage{amscd}
\usepackage{amsmath}
\usepackage{epsfig}
\usepackage{color}
\usepackage{ulem}

\def\be{\begin{equation}}
\def\ee{\end{equation}}
\def\ben{\begin{displaymath}}
\def\een{\end{displaymath}}
\def\baa{\begin{eqnarray}}
\def\eaa{\end{eqnarray}}

\def\ba{\begin{array}}
\def\ea{\end{array}}

\makeatletter
\@addtoreset{equation}{section}
\makeatother

\usepackage{amsthm}
\usepackage{amsmath}
\usepackage{amssymb}
\usepackage{amsbsy}
\usepackage{amsfonts}
\usepackage{latexsym}
\usepackage{amscd}
\usepackage{eucal}
\usepackage{mathrsfs}
\usepackage[all, knot]{xy}
\xyoption{arc} \xyoption{color}\xyoption{line}

\newtheorem{Remark}{Remark}

\def\be{\begin{equation}}
\def\ee{\end{equation}}
\def\ben{\begin{displaymath}}
\def\een{\end{displaymath}}
\def\baa{\begin{eqnarray}}
\def\eaa{\end{eqnarray}}

\def\ba{\begin{array}}
\def\ea{\end{array}}

\def\2x2{{\left(\!\!\begin{array}{cc}a&b\\c&d\\\end{array}\!\!\right)}}

\begin{document}
\title{Determinant of Laplacian on tori of constant positive curvature with one conical point}

\author{Victor Kalvin,  Alexey Kokotov}
\maketitle

\vskip2cm
{\bf Abstract.}
 We find an explicit expression for the zeta-regularized determinant  of (the Friedrichs extensions) of the Laplacians on a compact Riemann surface of genus one with conformal metric of curvature $1$ having a single conical singularity of angle $4\pi$.
\vskip2cm

\section{Introduction}

Let $X$ be a compact Riemann surface of genus one and let $P\in X$. According to \cite{CLW}, Cor. 3. 5. 1, there exists {\it at most} one conformal metric on $X$ of constant curvature $1$ with a (single) conical point of angle $4\pi$ at $P$. The following simple construction shows that such a metric, $m(X, P)$, in fact always exists (and due to \cite{CLW} is unique).

 Consider the spherical triangle $T=\{(x_1, x_2, x_3)\in S^2\subset {\mathbb R}^3: x_1\geq 0, x_2\geq 0, x_3\geq 0\}$ with  all three angles equal to $\pi/2$. Gluing two copies of $T$ along their boundaries, we get the Riemann sphere ${\mathbb P}$ with metric $m$ of curvature $1$ and three conical points $P_1, P_2, P_3$ of conical angle $\pi$. Consider the two-fold  covering
 \begin{equation}\mu: X(Q) \to {\mathbb P}   \end{equation}
ramified over $P_1$, $P_2$, $P_3$ and some point $Q\in {\mathbb P}\setminus \{P_1, P_2, P_3\}$.
Lifting the metric $m$ from ${\mathbb P}$ to the compact Riemann surface $X(Q)$ of genus one via $\mu$, one gets the metric $\mu^*m$ on $X(Q)$ which has curvature $1$ and the unique conical point of angle $4\pi$ at the preimage $\mu^{-1}(Q)$ of $Q$.
Clearly, any compact surface of genus one is (biholomorphically equivalent to) $X(Q)$ for some $Q\in {\mathbb P}\setminus \{P_1, P_2, P_3\}$.
Now let $X$ be an arbitrary compact Riemann surface of genus one and let $P$ be any point of $X$.
Take $Q\in {\mathbb P}$ such that $X=X(Q)$ and consider the automorphism $\alpha:X\to X$ (the translation) of $X$ sending $P$ to $\mu^{-1}(Q)$.
Then $$m(X, P)=\alpha^*(\mu^*(m))=(\mu\circ \alpha)^*(m)\,.$$

Now introduce the scalar (Friedrichs) self-adjoint Laplacian $\Delta(X, P):=\Delta^{m(X, P)}$ on $X$ corresponding to the metric $m(X, P)$. For any $P$ and $Q$ from $X$ the operators $\Delta(X, P)$ and $\Delta(X, Q)$ are isospectral and, therefore,
the $\zeta$-regularized (modified, i. e. with zero modes excluded) determinant ${\rm det}\Delta(X, P)$ is independent of $P\in X$ and, therefore, is a function on moduli space ${\cal M}_1$ of Riemann surfaces of genus one.
The main result of the present work is the following explicit formula for this function:
\begin{equation}\label{result}
{\rm det}\Delta(X, P) = C_1\,|\Im \sigma||\eta(\sigma)|^4F(t)=C_2\,{\rm det}\Delta^{(0)}(X)F(t),
\end{equation}
 where $\sigma$ is the $b$-period of the Riemann surface $X$, $C_1$ and $C_2$ are absolute constants,
 $\eta$ is the Dedekind eta-function, $\Delta^{(0)}$ is the Lapalacian on $X$ corresponding to the flat conformal metric of unit volume,
 the surface $X$ is represented as the two-fold covering of the Riemann sphere ${\mathbb C}P^1$ ramified over the poits $0, 1, \infty$ and $t\in {\mathbb C}\setminus \{0, 1\}$, and
 \begin{equation}
 F(t)=\frac{|t|^{\frac{1}{24}}|t-1|^{\frac{1}{24}}}{(|\sqrt{t}-1|+|\sqrt{t}+1|)^{\frac{1}{4}}}\,.
 \end{equation}

As it is well-known, the moduli space ${\cal M}_1$ coincides with the quotient space
$$\left( {\mathbb C}\setminus \{0, 1\}\right)/G\,,$$
where $G$ is a finite group of order $6$, generated by transformations $t\to \frac{1}{t}$ and $t\to 1-t$.
A direct 
check shows that
$F(t)=F(\frac{1}{t})$ and $F(t)=F(1-t)$ and, therefore, the right hand side of (\ref{result}) is in fact a function on ${\cal M}_1$.
\begin{Remark} Using the classical relation (see, e. g. \cite{Clemens} (3.35))
$$t=-\left(\frac{\Theta[^1_0](0\,|\,\sigma)}{\Theta[^0_1](0\,|\, \sigma)}\right)^4\,,$$
one can rewrite the right hand side as a function  $\sigma$ only.
\end{Remark}

The classical (see \cite{Polch}) relation ${\rm det}\Delta^{(0)}=C\,|\Im \sigma||\eta(\sigma)|^4$ used in (\ref{result}), implies that (\ref{result}) can be considered as a version of Polyakov's formula (relating determinants of the Laplacians corresponding to two {\it smooth} metrics in the same conformal class) for the case of two conformally equivalent metrics on a torus: one of them is smooth and flat, another is of curvature one and has one (very special) singular point.
\section{Metrics on the base and on the covering}
Here we find an explicit expression for the metric $m$ on the Riemann sphere ${\mathbb P}={\mathbb C}P^1$ of curvature $1$ and with three conical singularities at $P_1=0$, $P_2=1$ and $P_3=\infty$.

The stereographic projection (from the south pole) maps the spherical triangle $T$ onto quarter of the unit disk $\{z\in  {\mathbb C}; |z|\leq 1,\  0\leq {\rm Arg}\, z\leq \pi/2\}$. The conformal map
\begin{equation}\label{map}
z\mapsto w=\left(\frac{1+z^2}{1-z^2}\right)^2
\end{equation}
sends this quarter of the disk to the upper half-plane $H$; the corner points $i, 0, 1$ go to the points $0, 1$ and $\infty$ on the real line. The push forward of the standard round metric
$$\frac{4|dz|^2}{(1+|z|^2)^2}$$
on the sphere by this map gives rise to the metric
\begin{equation}\label{basemetric}m=\frac{|dw|^2}{|w||w-1|(|\sqrt{w}+1|+|\sqrt{w}-1|)^2}\end{equation}
on $H$; clearly, the latter metric can be extended (via the same formula) to ${\mathbb C}P^1$. The resulting curvature one metric on  ${\mathbb C}P^1$ (also denoted by $m$)  has three conical singularities of angle $\pi$: at $w=0$, $w=1$ and $ w=\infty$.

Consider a two-fold covering of the Riemann sphere by a compact Riemann surface $X(t)$ of genus $1$
\begin{equation}\label{cov}
\mu: X(t)\to {\mathbb C}P^1\end{equation}
ramified over four points: $0, 1, \infty$ and $t\in {\mathbb C}\setminus \{0, 1\}$.
Clearly, the pull back metric $\mu^*m$ on $X(t)$ is a curvature one metric with exactly one conical singularity. The singularity is a conical point of angle $4\pi$ located at  the point $\mu^{-1}(t)$.
\section{Determinant of Laplacian as  function of critical value  $t$}
The analysis from~\cite{KK} in particular implies that one can introduce the standard Ray-Singer $\zeta$-regularized determinant of the (Friedrichs) self-adjoint Laplacian
$\Delta^{\mu^*m}$ in $L_2(X(t), \mu^*m)$
$${\rm det}\,\Delta^{\mu^*m}:=\exp\{-\zeta'_{\Delta^{\mu^*m}}(0)\}\,,$$
where $\zeta'_{\Delta^{\mu^*m}}$ is the operator zeta-function.
In this section we establish a formula for the variation of  $\zeta'_{\Delta^{\mu^*m}}(0)$ with respect to the parameter $t$ (the fourth ramification point of the covering~\eqref{cov}. The derivation of this formula coincides almost verbatim with the proof of \cite[Proposition 6.1]{KK}, therefore, we will give only few details.

For the sake of brevity we identify the point $t$ of the base ${\mathbb C}P^1$ with its (unique) preimage $\mu^{-1}(t)$ on
$X(t)$.

Let $Y(\lambda; \,\cdot\, )$ be the (unique) special solution
of the Helmholz equation (here $\lambda$ is the spectral parameter)
$(\Delta^m -\lambda)Y=0$
on $X\setminus \{t\}$ with asymptotics $Y(\lambda) (x)=\frac{1}{x}+O(x)$ as $x\to 0$, where
$x(P)=\sqrt{\mu(P)-t}$ is {\ it the distinguished holomorphic local parameter} in a vicinity of the ramifiction point $t\in X(t)$ of the covering (\ref{cov}).
Introduce the complex-valued function $\lambda\mapsto b(\lambda)$ as the coefficient near $x$ in the asymptotic expansion
$$Y(x, \bar x; \lambda)=\frac{1}{x}+c(\lambda)+a(\lambda)\bar x+b(\lambda)x +O(|x|^{2-\epsilon})
\text {\ as } x\to 0.$$

The following variational formula is proved in \cite[Propositon 6.1]{KK}:
\begin{equation}\label{var1}
\partial_t (-\zeta'_{\Delta^{\mu^*m}}(0))=\frac{1}{2}\left(b(0)-b(-\infty)\right).
\end{equation}
The value $b(0)$ is found in \cite[Lemma 4.2]{KK}: one has the relation
\begin{equation}\label{me}b(0)=-\frac{1}{6}S_{Sch}(x)\Big|_{x=0}, \end{equation}
 where $S_{Sch}$ is the Schiffer projective connection on the Riemann surface $X(t)$.

 Since $\lambda=-\infty$ is a local regime, in order to find $b(-\infty)$ the solution $Y$ can be replaced by a local solution with the same asymptotic as $x\to 0$. A local solution $\widehat Y$ with asymptotic
 $$\widehat Y(u, \bar u; \lambda)=\frac{1}{u}+\hat c(\lambda)+\hat a(\lambda)\bar u+\hat b(\lambda)u +O(|u|^{2-\epsilon})
\text {\ as } u\to 0$$
 in the local parameter $u^2=z-s$ was constructed in~\cite[Lemma 4.1]{KK} by separation of variables; here $z$ and  $w=\mu(P)$ (resp. $s$ and $t$) are related by~\eqref{map} (resp. by~\eqref{map} with $z=s$ and $w=t$) and $\hat b(-\infty)=\frac{1}{2}\frac{\bar s}{1+|s|^2}$. One can easily find the coefficients $A(t)$ and $B(t)$ of the Taylor series  $u=A(t)x+B(t)x^3+O(x^5)$. As a local solution replacing $Y$ we can take $A(t)\widehat Y$. This immediately implies  $b(-\infty)=A^2(t)\hat b(-\infty)-B(t)/A(t)$. A straightforward  calculation verifies that
 \begin{equation}\label{Victor}
 b(-\infty)=\partial_t \log \left ({|t||t-1|(|\sqrt{t}+1|+|\sqrt{t}-1|)^2}\right )^{1/4}.
   \end{equation}
Observe that  the right hand side in~\eqref{Victor}  is actually the value of $\partial_w \log \rho(w, \bar w)^{-1/4}$  at $w=t$, where $\rho(w, \bar w)$ is the conformal factor of the metric  \eqref{basemetric}; this is also a direct consequence of \cite[Lemma 4]{VK}.

Using (\ref{var1}) together with (\ref{me}) and (\ref{Victor}), we are now able to derive an explicit formula for ${\rm det}\Delta^{\mu^*m}$.
\section{Explicit formula for the determinant}
Equations (\ref{var1}), (\ref{me}) and (\ref{Victor}) imply
that the determinant of the Laplacian ${\rm det}\,\Delta^{\mu^*m}=\exp\{-\zeta'_{\Delta^{\mu^*m}}(0)\}$
can be represented as a product
\begin{equation}\label{prelim}
{\rm det}\,\Delta^{\mu^*m}=C\,|\Im \sigma||\tau(t)|^2\left| \frac{1}{|t||t-1|(|\sqrt{t}+1|+|\sqrt{t}-1|)^2}\right|^{1/8}\,
\end{equation}
where $\tau(t)$ is the value of the Bergman tau-function (see \cite{KokKor1}, \cite{KokKor2}, \cite{KokStra}) on the Hurwitz
space $H_{1, 2}(2)$ of  two-fold genus one coverings of the Riemann sphere, having $\infty$ as a ramification point at the covering, ramified over $1, 0, \infty$ and $t$.
More specifically, $\tau$ is a solution of the equation
$$\frac{\partial \log \tau}{\partial t}=-\frac{1}{12}S_B(x)|_{x=0}\,$$
where $S_B$ is the Bergman projective connection on the covering Riemann surface $X(t)$ of genus one and $x$ is the distinguished holomorphic parameter in a vicinity of the ramification point $t$ of $X(t)$.
We remind the reader that
 the Bergman and the Schiffer projective connections are related via the equation
 $$S_{Sch}(x)=S_B(x)-6\pi(\Im\sigma)^{-1} v^2(x)\,$$
where $v$ is the normalized holomorphic differential on $X(t)$ and that the Rauch variational formula (see, e. g., \cite{KokKor1})
implies the relation
$$\frac {\partial\log \Im \sigma}{\partial t}=\frac{\pi}{2}(\Im \sigma)^{-1} v^2(x)|_{x=0}\,.$$

The needed explicit expression for $\tau$ can be found e. g. in \cite[ f-la (18)]{KokStra} (it is a very special case of the explicit formula for the Bergman tau-function on general coverings of arbitrary genus and degree found in \cite{KokKor2} as well as of a much earlier  formula of Kitaev and Korotkin for hyperelliptic coverings~\cite{KitKor}). Namely,  \cite[ f-la (18)]{KokStra}  implies that
\begin{equation}\label{tau}
\tau=\eta^2(\sigma)\left[\frac{v(\infty)^3}{v(P_1)v(P_2)v(Q)}\right]^{\frac{1}{12}},
\end{equation}
where $P_1$ and $P_2$ are the points of the $X(t)$ lying over $0$ and $1$, $Q$ is a point of $X(t)$ lying over $t$
 and  $\infty$ denotes the point of the covering curve $X(t)$ lying over the point at infinity of the base ${\mathbb C}P^1$;
 $v$ is an arbitrary nonzero holomorphic differential on $X(t)$; and, say,  $v(P_1)$ is the value of this differential in the distinguished holomorphic parameter at $P_1$. (One has to  take into account that  $\tau=\tau_I^{-2}$, where $\tau_I$ is from \cite{KokStra}.)
Taking
$$v=\frac{dw}{\sqrt{(w(w-1)(w-t)}}\,,$$
and using the following expressions for the distinguished local parameters at $P_1$, $P_2$, $Q$ and $\infty$
$$x=\sqrt{w}; \ \ x=\sqrt{w-1};\ \ x=\sqrt{w-t}; \ \ x=\frac{1}{\sqrt w}$$
one arrives at the relations (where $\sim$ means $=$ up to insignificant constants like $\pm 2$, etc.)
$$v(P_1)\sim\frac{1}{\sqrt{t}}; \ \ v(P_2)\sim\frac{1}{\sqrt{t-1}}; \ \  v(Q)\sim\frac{1}{\sqrt{t(t-1)}}; \ \ v(\infty)\sim1.$$

These relations together with (\ref{tau}) and (\ref{prelim}) imply (\ref{result}).


\begin{thebibliography}{100}

\bibitem{CLW} Ching-Li Chai, Chang-Shou Lin, Chin-Lung Wang, Mean field equation, hyperelliptic curves and modular forms: I,
Cambridge Journal of Mathematics, Vol. 3, N 1-2, 2015
\bibitem{Clemens} C. Clemens, A scrapbook of complex curve theory, Grad. Studies in Math., Vol 55
\bibitem{KK} V. Kalvin, A. Kokotov, Metrics of constant positive curvature, Hurwitz spaces and ${\rm det \Delta}$, IMRN, 2018; in press; arXiv:1612.08660
\bibitem{KokKor1} A. Kokotov, D. Korotkin, Tau-functions on Hurwitz spaces,
Mathematical Physics,
Analysis and Geometry, 7 (2004), no. 1, 47--96.
\bibitem{KokKor2} A. Kokotov, D. Korotkin, Isomonodromic tau-function of Hurwitz Frobenius manofolds,  Int. Math. Res. Not. IMRN (2006), pp. 1-34
\bibitem{KokStra} A.Kokotov, I. Strachan,  On the isomonodromic tau-function for the Hurwitz spaces of branched coverings of genus zero and one, Mathematical Research Letters, 12, 2005,  no. 5-6, 857–-875.

\bibitem{KitKor} V. Kitaev, D. Korotkin,
On solutions of the Schlesinger equations in terms of
theta-functions,
International Mathematics Research Notices,
1998,
no. 17,877–-905.
\bibitem{VK} V. Kalvin, On Determinants of Laplacians on Compact Riemann Surfaces Equipped with Pullbacks of  Conical Metrics by Meromorphic Functions, in preparation
\bibitem{Polch}  J. Polchinski,  Evaluation of the one loop string path integral. Comm. Math. Phys. 104 (1986), no. 1, 37--47

\end{thebibliography}
\end{document}